\documentclass[12pt,fleqn,leqno]{amsart}
\usepackage{amsfonts,amssymb}
\usepackage{mathrsfs}
\usepackage{diagram}

\numberwithin{equation}{section}

\newtheorem{theorem}{Theorem}[section]
\newtheorem{lemma}[theorem]{Lemma}
\newtheorem{corollary}[theorem]{Corollary}
\newtheorem{proposition}[theorem]{Proposition}

\renewcommand{\emptyset}{\varnothing}
\DeclareMathOperator{\graph}{Graph}
\DeclareMathOperator{\var}{var}
\DeclareMathOperator{\uhr}{\upharpoonright}
\DeclareMathOperator{\diam}{diam}
\DeclareMathOperator{\mesh}{mesh}
\DeclareMathOperator{\hgt}{ht}
\DeclareMathOperator{\td}{td}
\DeclareMathOperator{\I}{\mathbb I}
\DeclareMathOperator{\dom}{dom}
\DeclareMathOperator{\R}{\mathbb{R}}
\DeclareMathOperator{\sphere}{\mathbb{S}}
\DeclareMathOperator{\B}{\mathbb{B}}
\DeclareMathOperator{\smap}{\rightsquigarrow}

\begin{document}


\title[Open maps having the Bula property] {Open maps having the Bula
  property}

\author{Valentin Gutev}

\address{School of Mathematical Sciences, University of KwaZulu-Natal,
  Westville Campus, Private Bag X54001, Durban 4000, South Africa}

\email{gutev@ukzn.ac.za}

\thanks{Research of the first author is supported in part by the NRF
  of South Africa.}

\author{Vesko Valov}

\address{Department of Computer Science and Mathematics, Nipissing
  University, 100 College Drive, P.O. Box 5002, North Bay, ON, P1B
  8L7, Canada}

\email{veskov@nipissingu.ca}

\thanks{The second author was partially supported by NSERC Grant
  261914-03.}

\keywords{Dimension, $C$-space, set-valued mapping, selection}

\subjclass[2000]{Primary 54F45, 54F35, 54C60, 54C65; Secondary 55M10,
  54C35, 54B20}


\begin{abstract}
  Every open continuous map $f$ from a space $X$ onto a paracompact
  $C$-space $Y$ admits two disjoint closed sets $F_0,F_1\subset X$, with
  $f(F_0)=Y=f(F_1)$, provided all fibers of $f$ are infinite and
  $C^*$-embedded in $X$. Applications are demonstrated for the
  existence of ``disjoint'' usco multiselections of set-valued l.s.c.\
  mappings defined on paracompact $C$-spaces, and for special type of
  factorizations of open continuous maps from metrizable spaces onto
  paracompact $C$-spaces. This settles several open questions raised
  in \cite{karassev-tuncali-valov:07}.
\end{abstract}
\date{\today}
\maketitle


\section{Introduction}
\label{section-introduction}

All spaces in this paper are assumed to be at least completely
regular. Following Kato and Levin \cite{kato-levin:00}, a continuous
surjective map $f\colon X\to Y$ is said to have the \emph{Bula property} if
there exist two disjoint closed subsets $F_0$ and $F_1$ of $X$ such
that $f(F_0)=Y=f(F_1)$. In the sequel, such a pair $(F_0,F_1)$ will be
called a \emph{Bula pair} for $f$. Bula \cite{bula:83} proved that
every open continuous map $f$ from a compact Hausdorff space onto a
finite-dimensional metrizable space has this property provided all
fibers of $f$ are dense in themselves. This result was generalized in
\cite{gutev:93b} to the case $Y$ is countable-dimensional. Recently,
Levin and Rogers \cite{levin-rogers:99} obtained a further
generalization with $Y$ being a $C$-space. The question whether the
Levin-Rogers result remains true for open maps between metrizable
spaces was raised in \cite[Problem 1514]{karassev-tuncali-valov:07}
(if $Y$ is not a $C$-space, this is not true, see
\cite{dranishnikov:90} and \cite{levin-rogers:99}). Here, we provide a
positive answer to this question:

\begin{theorem}
  \label{theorem-bula-property-basic}
  Let $X$ be a space, $Y$ be a paracompact $C$-space, and let $f\colon X\to
  Y$ be an open continuous surjection such that all fibers of $f$ are
  infinite and $C^*$-embedded in $X$. Then, $f$ has the Bula property.
\end{theorem}

The $C$-space property was originally defined by W.\ Haver
\cite{haver:1974} for compact metric spaces. Later on, Addis and
Gresham \cite{addis-gresham:78} reformulated Haver's definition for
arbitrary spaces: A space $X$ has property $C$ (or $X$ is a
\emph{$C$-space}) if for every sequence $\{\mathscr{W}_n:n=1,2,\dots\}$
of open covers of $X$ there exists a sequence
${\{\mathscr{V}_n:n=1,2,\dots\}}$ of pairwise disjoint open families in
$X$ such that each $\mathscr{V}_n$ refines $\mathscr{W}_n$ and
$\bigcup\{\mathscr{V}_n:n=1,2,\dots\}$ is a cover of $X$. It is well-known
that every finite-dimensional paracompact space, as well as every
countable-dimensional metrizable space, is a $C$-space
\cite{addis-gresham:78}, but there exists a compact metric $C$-space
which is not countable-dimensional \cite{pol:81}.  Let us also remark
that a $C$-space $X$ is paracompact if and only if it is countably
paracompact and normal. Finally, let us recall that a subset $A\subset X$ is
\emph{$C^*$-embedded} in $X$ if every bounded real-valued continuous
function on $A$ is continuously extendable to the whole of
$X$. \medskip

Theorem \ref{theorem-bula-property-basic} has several interesting
applications. In Section \ref{section-bula-pairs0selections}, we apply
this theorem to the graph of l.s.c.\ set-valued mappings defined on
paracompact $C$-spaces and with point-images being closed and infinite
subsets of completely metrizable spaces. Thus, we get that any such
l.s.c.\ mapping has a pair of ``disjoint'' usco multiselections (see,
Corollaries \ref{corollary-strict-michael-pair} and
\ref{corollary-disjoint-michael-pair}), which provides the complete
affirmative solution to \cite[Problem 1515]{karassev-tuncali-valov:07}
and sheds some light on \cite[Problem
1516]{karassev-tuncali-valov:07}. In this regard, let us stress the
attention that, in Theorem \ref{theorem-bula-property-basic}, no
restrictions on $X$ are called a priori. In Section
\ref{section-open-maps-projections}, we use Theorem
\ref{theorem-bula-property-basic} to demonstrate that every open
continuous map $f$ from a metric space $(X,d)$ onto a paracompact
$C$-space $Y$ admits a special type of factorization ($Y\times[0,1]$
throughout), provided all fibers of $f$ are dense in themselves and
complete with respect to $d$, see Theorem
\ref{theorem-cantor-bouquet}. This result is a common generalization
of \cite[Theorem 1.1]{gutev:93b} and \cite[Theorem
1.2]{levin-rogers:99} (see, Corollary
\ref{corollary-closed-section-zero-compact}), and provides the
complete affirmative solution to \cite[Problem
1512]{karassev-tuncali-valov:07}.  \medskip

Finally, a word should be said also for the proof of Theorem
\ref{theorem-bula-property-basic} itself. Briefly, a preparation for
this is done in the next section. It is based on the existence of a
continuous function $g:X\to [0,1]$ such that $g$ is not constant on each
fiber $f^{-1}(y)$, $y\in Y$, of $f$ (see, Theorem
\ref{theorem-fiber-constant}).  Having already established this, the
proof of Theorem \ref{theorem-bula-property-basic} will be
accomplished in Section \ref{section-proof-main-theorem} relying on a
``parametric'' version of an idea in the proof of \cite[Theorem
1.3]{levin-rogers:99}.

\section{Bula property and fiber-constant maps}
\label{section-bula-property}

Suppose that $(F_0,F_1)$ is a Bula pair for a map $f:X\to Y$, where $X$
is a normal space. Then, there exists a continuous function $g:X\to
[0,1]$ such that $g\uhr f^{-1}(y)$ is not constant for every $y\in
Y$. Indeed, take $g:X\to [0,1]$ to be such that $F_i\subset g^{-1}(i)$,
$i=0,1$.  In this section, we demonstrate that the map $f$ in Theorem
\ref{theorem-bula-property-basic} has this property as well. Namely,
the following theorem will be proved.

\begin{theorem}
  \label{theorem-fiber-constant}
  Let $X$ be a space, $Y$ be a paracompact $C$-space, and let $f\colon X\to
  Y$ be an open continuous surjection such that all fibers of $f$ are
  infinite and $C^*$-embedded in $X$. Then, there exists a continuous
  function $g:X\to [0,1]$ such that $g\uhr f^{-1}(y)$ is not constant
  for every $y\in Y$.
\end{theorem}

To prepare for the proof of Theorem \ref{theorem-fiber-constant}, let
us recall some terminology. For spaces $Y$ and $Z$, we will use
$\Phi:Y\smap Z$ to denote that $\Phi$ is a set-valued mapping, i.e.\ a map
from $Y$ into the nonempty subsets of $Z$. A mapping $\Phi:Y\smap Z$ is
\emph{lower semi-continuous}, or l.s.c., if the set
\[
\Phi^{-1}(U)=\{y\in Y:\Phi(y)\cap U\neq\emptyset\}
\]
is open in $Y$ for every open $U\subset Z$.  A mapping $\Phi:Y\smap Z$ has an
\emph{open} (\emph{closed}) graph if its graph
\[
\graph(\Phi)=\big\{(y,z)\in Y\times Z: z\in \Phi(y)\big\}
\]
is open (respectively, closed) in $Y\times Z$. A map $g:Y\to Z$ is a
\emph{selection} for $\Phi:Y\smap Z$ if $g(y)\in \Phi(y)$ for every $y\in
Y$. Finally, let us recall that a space $Z$ is $C^m$ for some
$m\geq 0$ if every continuous image of the $k$-dimensional sphere
$\sphere^{k}$ ($k\leq m$) in $Z$ is contractible in $Z$.\medskip

In what follows, $\I=[0,1]$ and $C(Z,\I)$ denotes the set of all
continuous functions from $Z$ to $\I$. Also, $C(Z)=C(Z,\R)$ is the
set of all continuous functions on $Z$, and $C^*(Z)$ --- that of all
bounded members of $C(Z)$. As usual, $C^*(Z)$ is equipped  with the
\emph{sup-metric} $d$ defined by
\[
d(g,h)=\sup\big\{|g(z)-h(z)| :z\in Z\big\},\quad g,h\in C^*(Z).
\]
It should be mentioned that $C^*(Z)$ is a Banach space, and
$C(Z,\I)$ is a closed convex subset of $C^*(Z)$. In the sequel, for
$g\in C(Z,\I)$ and $\varepsilon>0$, we will use $B_\varepsilon^d(g)$
to denote the open $\varepsilon$-ball $B_\varepsilon^d(g)=\{h\in
C(Z,\I): d(g,h)<\varepsilon\}$.
\medskip

The next statement is well-known and easy to prove.

\begin{lemma}
  \label{proposition-restriction-map-open}
    Let $X$ be a space, and let $A\subset X$ be a $C^*$-embedded
  subset of $X$. Then, the restriction map $\pi_A:C(X,\I)\to C(A,\I)$ is
  an open continuous surjection.
\end{lemma}

For a subset $B$ of a space $Z$, let $\Theta_Z(B,\I)$ be the set of all
members of $C(Z,\I)$ which are constant on $B$. If $B=Z$, then we will
denote this set merely by $\Theta(Z,\I)$. Note that $\Theta(Z,\I)$ is, in fact,
homeomorphic to $\I$.

\begin{proposition}
  \label{proposition-finite-contarctible-values}
  Let $X$ be a space, and let $A\subset X$ be an infinite $C^*$-embedded
  subset of $X$. Then, the set $C(X,\I)\setminus \Theta_X(A,\I)$ is $C^m$ for every
  $m\geq 0$.
\end{proposition}

\begin{proof}
  Consider the restriction map $\pi_A:C(X,\I)\to C(A,\I)$, and take a
  continuous map $g:\sphere^n\to C(X,\I)\setminus \Theta_X(A,\I)$ for some $n\geq
  0$. Then, by Lemma \ref{proposition-restriction-map-open}, the
  composition $\pi_A\circ g:\sphere^n\to C(A,\I)\setminus \Theta(A,\I)$ is also
  continuous. Observe that $C(A,\I)$ is an infinite-dimensional closed
  convex subset of $C^*(A)$ because $A$ is infinite.  From another
  hand, $\Theta(A,\I)$ is one-dimensional being homeomorphic to $\I$. Then,
  by \cite[Lemma 2.1]{michael:88a}, $C(A,\I)\setminus \Theta(A,\I)$ is $C^m$ for
  all $m\geq 0$. Hence, there exists a continuous extension $\ell:\B^{n+1}\to
  C(A,\I)\setminus \Theta(A,\I)$ of $\pi_A\circ g$ over the $(n+1)$-dimensional ball
  $\B^{n+1}$. Consider the set-valued mapping ${\Phi:\B^{n+1}\smap
    C(X,\I)}$ defined by $\Phi(t)=\{g(t)\}$ if $t\in \sphere^n$ and
  $\Phi(t)=\pi_A^{-1}(\ell(t))$ otherwise. Since $g$ is a selection for
  $\pi_A^{-1}\circ\ell \uhr \sphere^n$ and, by Lemma
  \ref{proposition-restriction-map-open}, the restriction map is $\pi_A$
  is open, the mapping $\Phi$ is l.s.c. (see, \cite[Examples 1.1$^*$ and
  1.3$^*$]{michael:56a}). Also, $\Phi$ is closed and convex-valued in
  $C(X,\I)$, hence in the Banach space $C^*(X)$ as well.  Then, by the
  Michael's selection theorem \cite[Theorem 3.2$''$]{michael:56a}, $\Phi$
  has a continuous selection ${h:\B^{n+1}\to C(X,\I)}$ which is, in
  fact, a continuous extension of $g$ over $\B^{n+1}$. Moreover,
  $\pi_A(h(t))=\ell(t)\notin \Theta(A,\I)$ for all $t\in \B^{n+1}$, which completes the
  proof.
\end{proof}

A function $\xi:X\to \R$ is \emph{lower} (\emph{upper})
\emph{semi-continuous} if the set
\[
\{x\in x:\xi(x)>r\}\quad \text{(respectively,
  $\{x\in X:\xi(x)<r\}$)}
\]
is open in $X$ for every $r\in \R$. Suppose that $f:X\to Y$ is a
surjective map. Then, to any $g:X\to \I$ we will associate the
functions ${\inf[g,f],\sup[g,f]:Y\to \I}$ defined for $y\in Y$ by
\[
\inf[g,f](y)=\inf\big\{g(x): x\in f^{-1}(y)\big\},
\]
and, respectively,
\[
\sup[g,f](y)=\sup\big\{g(x): x\in f^{-1}(y)\big\}.
\]
Finally, we will also associate the function $\var[g,f]:X\to \I$ defined
by
\[
\var[g,f](y)=\sup[g,f](y)-\inf[g,f](y),\quad  y\in Y.
\]
Observe that $g:X\to \I$ is not constant on any fiber $f^{-1}(y)$, $y\in
Y$, if and only if $\var[g,f]$ is positive-valued. The following
property is well-known, \cite{isiwata:67} (see, also,
\cite[1.7.16]{engelking:89}).

\begin{proposition}
  [\cite{isiwata:67}]
  \label{proposition-lower-uper-variation}
  Let $X$ and $Y$ be spaces, $f:X\to Y$ be an open surjective map, and
  let $g\in C(X,\I)$. Then, $\sup[g,f]$ is lower semi-continuous, while
  $\inf[g,f]$ is upper semi-continuous. In particular, $\var[g,f]$ is
  lower semi-continuous.
\end{proposition}

We finalize the preparation for the proof of Theorem
\ref{theorem-fiber-constant} with the following proposition.

\begin{proposition}
  \label{proposition-closed-graph-variation}
  Let $X$ and $Y$ be spaces, and let $f:X\to Y$ be an open surjective
  map. Then the set-valued mapping $\Theta:Y\smap C(X,\I)$ defined by
  $\Theta(y)=\Theta_X(f^{-1}(y),\I)$, $y\in Y$, has a closed
  graph.
\end{proposition}

\begin{proof}
  Take a point $y\in Y$ and $g\notin \Theta(y)$. Then, $\var[g,f](y)> 2\delta$ for some
  positive number $\delta>0$. By Proposition
  \ref{proposition-lower-uper-variation}, there exists a neighbourhood
  $V$ of $y$ such that $\var[g,f](z)>2\delta$ for every $z\in V$. Then, $V\times
  B_\delta^d(g)$ is an open set in $Y\times C(X,\I)$ such that $\big(V\times
  B_\delta^d(g)\big)\cap \graph(\Theta)=\emptyset$. Indeed, take $z\in V$ and $h\in
  B_\delta^d(g)$. Since $\var[g,f](z)>2\delta$, there are points $x,t\in
  f^{-1}(z)$ such that $|g(x)-g(t)|>2\delta$. Since $h\in B_\delta^d(g)$, we have
  $|h(x)-g(x)|<\delta$ and $|h(t)-g(t)|<\delta$. Hence, $h(x)\neq h(t)$, which
  implies that $\var[h,f](z)>0$. Consequently, $h\notin \Theta(z)$.
\end{proof}

\begin{proof}[Proof of Theorem \ref{theorem-fiber-constant}]
  Consider the set-valued mapping $\Phi:Y\smap C(X,\I)$ defined by
  $\Phi(y)=C(X,\I)\setminus \Theta(y)$, $y\in Y$, where $\Theta$ is as in Proposition
  \ref{proposition-closed-graph-variation}. Then, by Proposition
  \ref{proposition-closed-graph-variation}, $\Phi$ has an open graph,
  while, by Proposition \ref{proposition-finite-contarctible-values},
  each $\Phi(y)$, $y\in Y$, is $C^m$ for all $m\geq 0$. Since $Y$ is a
  paracompact $C$-space, by the Uspenskij's selection theorem
  \cite[Theorem 1.3]{uspenskij:98}, $\Phi$ has a continuous selection
  $\varphi:Y\to C(X,\I)$. Define a map $g:X\to \I$ by $g(x)=[\varphi(f(x))](x)$, $x\in
  X$. Since $f$ and $\varphi$ are continuous, so is $g$ (see, the proof of
  \cite[Theorem 6.1]{gutev-valov:01a}). Since $g\uhr f^{-1}(y)=
  \varphi(y)\uhr f^{-1}(y)$ and $\varphi(y)\notin \Theta(y)$ for every $y\in Y$, $g$ is as
  required.
\end{proof}

\section{Proof of Theorem \ref{theorem-bula-property-basic}}
\label{section-proof-main-theorem}

Suppose that $X$, $Y$ and $f:X\to Y$ are as in Theorem
\ref{theorem-bula-property-basic}. By Theorem
\ref{theorem-fiber-constant}, there exists a function $g\in C(X,\I)$
such that $\inf[g,f](y)<\sup[g,f](y)$ for every $y\in Y$. Since
$\inf[g,f]$ is upper semi-continuous and $\sup[g,f]$ is lower
semi-continuous (by Proposition
\ref{proposition-lower-uper-variation}), and $Y$ is paracompact, by a
result of \cite{dieudonne:44} (see, also,
\cite{dowker:51,katetov:51}), there are continuous functions
$\gamma_0,\gamma_1:Y\to \I$ such that
  \[
  \inf[g,f](y)<\gamma_0(y)<\gamma_1(y)<\sup[g,f](y),\quad y\in Y.
  \]
  Let $\alpha_i=\gamma_i\circ f:X\to \I$, $i=0,1$. Then,
  \begin{equation}
    \label{eq:separating-variation}
    \inf[g,f](f(x))<\alpha_0(x)<\alpha_1(x)<\sup[g,f](f(x))\quad \text{for every $x\in X$.}
  \end{equation}
  Next, define a continuous function $\ell:X\times \I\to \R$ by letting 
  \[
  \ell(x,t)=\frac{t-\alpha_0(x)}{\alpha_1(x)-\alpha_0(x)},\quad (x,t)\in X\times \I.
  \]
  Observe that $\ell(x,\alpha_0(x))=0$ and $\ell (x,\alpha_1(x))=1$ for every $x\in
  X$. Hence,
  \begin{equation}
    \label{eq:linear-inequality}
    \ell\big(\{x\}\times [\alpha_0(x),\alpha_1(x)]\big)= [0,1],\quad \text{for every $x\in X$,}
  \end{equation}
  because $\ell$ is linear for every fixed $x\in X$. Finally, define a
  continuous function $h:X\to \R$ by $h(x)=\ell(x,g(x))$, $x\in X$. According
  to \eqref{eq:separating-variation} and \eqref{eq:linear-inequality},
  we now have that, for every $y\in Y$,
  \[
  h^{-1}\big((-\infty,0]\big)\cap f^{-1}(y)\neq\emptyset\neq h^{-1}\big([1,+\infty)\big)\cap
  f^{-1}(y).
  \]
  Then, $F_0=h^{-1}\big((-\infty,0]\big)$ and $F_1=h^{-1}\big([1,+\infty)\big)$
  are as required. The proof of Theorem
  \ref{theorem-bula-property-basic} completes.

\section{Bula pairs  and multiselections}
\label{section-bula-pairs0selections}

A set-valued mapping $\varphi:Y\smap Z$ is called a \emph{multiselection}
for $\Phi:Y\smap Z$ if $\varphi(y)\subset \Phi(y)$ for every $y\in Y$. In this section, we
present several applications of Theorem
\ref{theorem-bula-property-basic} about multiselections of l.s.c.\
mappings based on the following consequence of it.

\begin{corollary}
  \label{corollary-disjoint-lsc}
  Let $Y$ be a paracompact $C$-space, $Z$ be a normal space, and let
  $\Phi:Y\smap Z$ be an l.s.c.\ mapping such that each $\Phi(y)$, $y\in Y$, is
  infinite and closed in $Z$. Then, there exists a closed-graph
  mapping $\theta:Y\smap Z$ such that $\Phi(y)\cap \theta(y)\neq\emptyset\neq \Phi(y)\setminus \theta(y)$ for
  every $y\in Y$. 
\end{corollary}

\begin{proof}
  Let $X=\graph(\Phi)$ be the graph of $\Phi$, and let $f:X\to Y$ be the
  projection. Then, $f$ is an open continuous map (because $\Phi$ is
  l.s.c.) such that all fibers of $f$ are infinite. Let us observe
  that each $f^{-1}(y)$, $y\in Y$, is $C^*$-embedded in $X$. Indeed,
  take a point $y\in Y$, and a continuous function $g:f^{-1}(y)\to
  \I$. Since $f^{-1}(y)=\{y\}\times \Phi(y)$, we may consider the continuous
  function $g_0:\Phi(y)\to \I$ defined by $g_0(z)=g(y,z)$, $z\in \Phi(y)$. Since
  $Z$ is normal, there exists a continuous extension $h_0:Z\to \I$ of
  $g_0$. Finally, define $h:X\to \I$ by $h(t,z)=h_0(z)$ for every $t\in Y$
  and $z\in \Phi(t)$. Then, $h$ is a continuous extension of $g$. Thus, by
  Theorem \ref{theorem-bula-property-basic}, there are disjoint closed
  subsets $F_0, F_1\subset X$ such that $f(F_0)=Y=f(F_1)$. Finally, take a
  closed set $F\subset Y\times Z$, with $F\cap X=F_0$, and define $\theta:Y\smap Z$ by
  $\graph(\theta)=F$. This $\theta$ is as required.
\end{proof}

To prepare for our applications, we need also the following
observation about l.s.c.\ multiselections of l.s.c.\ mappings.

\begin{proposition}
  \label{proposition-lsc-selection-closed}
  Let $Y$ be a paracompact space, $Z$ be a space, $\Phi:Y\smap Z$
  be an l.s.c.\ closed-valued mapping, and let $\Psi:Y\smap Z$ be an
  open-graph mapping, with $\Phi(y)\cap \Psi(y)\neq\emptyset$ for every $y\in Y$.  Then,
  there exists a closed-valued l.s.c.\ mapping $\varphi:Y\smap Z$ such that
  $\varphi(y)\subset \Phi(y)\cap \Psi(y)$ for every $y\in Y$.
\end{proposition}

\begin{proof}
  Whenever $y\in Y$, there are open sets $V_y\subset Y$ and $W_y\subset Z$ such that
  $y\in V_y\subset \Phi^{-1}(W_y)$ and $V_y\times\overline{W_y}\subset \graph(\Psi)$. Indeed,
  take a point $z\in \Phi(y)\cap \Psi(y)$. Since $\Psi$ has an open graph, there are
  open sets $O_y\subset Y$ and $W_y\subset Z$ such that $y\in O_y$, $z\in W_y$ and
  $O_y\times\overline{W_y}\subset \graph(\Psi)$. Then, $V_y=O_y\cap \Phi^{-1}(W_y)$ is as
  required. Now, for every $y\in Y$, define a closed-valued mapping
  $\varphi_y:V_y\smap Z$ by letting that $\varphi_y(t)=\overline{\Phi(t)\cap W_y}$, $t\in
  V_y$. According to \cite[Propositions 2.3 and 2.4]{michael:56a},
  each $\varphi_y$, $y\in Y$, is l.s.c. Next, using that $Y$ is paracompact,
  take a locally-finite open cover $\mathscr{U}$ of $Y$ refining
  $\{V_y:y\in Y\}$ and a map $p:\mathscr{U}\to Y$ such that $U\subset V_{p(U)}$,
  $U\in \mathscr{U}$. Finally, define a mapping $\varphi:Y\smap Z$ by letting
  that
  \[
  \varphi(y)=\bigcup\big\{\varphi_{p(U)}(y): U\in \mathscr{U}\ \text{and}\ y\in U\big\},\quad
  y\in Y.
  \]
  This $\varphi$ is as required.
\end{proof}

In what follows, a mapping $\psi :Y\smap Z$ is \emph{upper
  semi-continuous}, or u.s.c., if the set
\[
\Phi^{\#}(U)=\{y\in Y:\Phi(y)\subset U\}
\]
is open in $Y$ for every open $U\subset Z$. Motivated by \cite{michael:59},
we say that a pair $(\varphi,\psi)$ of set-valued mapping $\varphi,\psi:Y\smap Z$ is a
\emph{Michael pair} for $\Phi:Y\smap Z$ if $\varphi$ is compact-valued and
l.s.c., $\psi$ is compact-valued and u.s.c., and $\varphi(y)\subset \psi(y)\subset \Phi(y)$ for
every $y\in Y$. \medskip

The following consequence provides the complete affirmative solution
to \cite[Problem 1515]{karassev-tuncali-valov:07}.

\begin{corollary}
  \label{corollary-strict-michael-pair}
  Let $(Z,\rho)$ be metric space, $Y$ be a paracompact $C$-space, and let
  $\Phi:Y\smap Z$ be an l.s.c.\ mapping such that each $\Phi(y)$, $y\in Y$, is
  infinite and $\rho$-complete. Then $\Phi$ has Michael pair $(\varphi,\psi):Y\smap
  Z$ such that $\Phi(y)\setminus \psi(y)\neq \emptyset$ for every $y\in Y$.
\end{corollary}

\begin{proof}
  By Corollary \ref{corollary-disjoint-lsc}, there is a closed-graph
  mapping $\theta:Y\smap Z$ such that $\Phi(y)\cap \theta(y)\neq\emptyset\neq \Phi(y)\setminus \theta(y)$ for every
  $y\in Y$. Consider the set-valued mapping $\Psi:Y\smap Z$ defined by
  $\graph(\Psi)=(Y\times Z)\setminus \graph(\theta)$. On one hand, by the properties of
  $\theta$, we have that $\Phi(y)\cap \Psi(y)\neq \emptyset$ for every $y\in Y$. On another hand,
  $\Phi$ is closed-valued having $\rho$-complete values. Hence, by
  Proposition \ref{proposition-lsc-selection-closed}, there exists a
  closed-valued l.s.c.\ mapping $\Phi_0:Y\smap Z$ such that $\Phi_0(y)\subset
  \Phi(y)\cap \Psi(y)$ for every $y\in Y$. Then, $\Phi_0$ has also $\rho$-complete
  values and, by a result of \cite{michael:59}, it has a Michael pair
  $(\varphi,\psi)$. This $(\varphi,\psi)$ is as required.
\end{proof}

We conclude this section with the following further application of
Theorem \ref{theorem-bula-property-basic} that sheds some light on
\cite[Problem 1516]{karassev-tuncali-valov:07}.

\begin{corollary}
  \label{corollary-disjoint-michael-pair}
  Let $(Z,\rho)$ be metric space, $Y$ be a paracompact $C$-space, and
  let $\Phi:Y\smap Z$ be an l.s.c.\ mapping such that each $\Phi(y)$,
  $y\in Y$, is infinite and $\rho$-complete. Then $\Phi$ has Michael
  pairs $(\varphi_i,\psi_i):Y\smap Z$, $i=0,1$, such that
  $\psi_0(y)\cap \psi_1(y)=\emptyset$ for every $y\in Y$.
\end{corollary}

\begin{proof}
  According to Corollary \ref{corollary-strict-michael-pair}, $\Phi$ has
  a Michael pair $(\varphi_0,\psi_0):Y\smap Z$ such that $\Phi(y)\setminus \psi_0(y)\neq\emptyset$ for
  every $y\in Y$. Note that $\psi_0$ has a closed-graph being u.s.c. Then,
  just like in the proof of Corollary
  \ref{corollary-strict-michael-pair}, there exists a Michael pair
  $(\varphi_1,\psi_1):Y\smap Z$ for $\Phi$ such that $\psi_1(y)\subset \Phi(y)\setminus \psi_0(y)$, $y\in
  Y$. These $(\varphi_i,\psi_i)$, $i=0,1$, are as required.
\end{proof}

\section{Open maps looking like projections}
\label{section-open-maps-projections}

Throughout this section, by a \emph{dimension} of a space $Z$ we mean
the covering dimension $\dim(Z)$ of $Z$.  In particular, $Z$ is
\emph{$0$-dimensional} if $\dim(Z)=0$.\medskip

We say that a continuous map $f:X\to Y$ has \emph{dimension} $\leq k$ if
all fibers of $f$ have dimension $\leq k$. A continuous map $f:X\to Y$ is
\emph{light} if it is $0$-dimensional, i.e.\ if $f$ has
$0$-dimensional fibers. Also, for convenience, we shall say that a
map $f:X\to Y$ is \emph{compact} if each fiber $f^{-1}(y)$, $y\in Y$, is a
compact subset of $X$.\medskip

Suppose that $f:X\to Y$ is a surjective map. A subset $F\subset X$ will be
called a \emph{section} for $f$ if $f(F)=Y$. In particular, we shall
say that a section $F$ for $f$ is \emph{open} (\emph{closed}) if $F$
is an open (respective, a closed) subset of $X$.\medskip

 In this section, we demonstrate the following factorization theorem
 which is a partial generalization of \cite[Theorem 1.1]{gutev:93b},
 also it provides the complete affirmative solution to \cite[Problem
 1512]{karassev-tuncali-valov:07}.

 \begin{theorem}
   \label{theorem-cantor-bouquet}
   Let $(X,d)$ be a metric space, $Y$ be a paracompact $C$-space, and
   let $f:X\to Y$ be an open continuous surjection such that each fiber
   of $f$ is dense in itself and $d$-complete. Also, let $U\subset X$ be an
   open section for $f$. Then, there exists a continuous surjective
   map $g:X\to Y\times \I$, a closed section $H\subset X$ for $f$, with $H\subset U$, and
   a copy $\mathfrak{C}\subset \I$ of the Cantor set such that
   \begin{enumerate}
   \item[(a)] $f=P_Y\circ g$, where $P_Y:Y\times \I\to Y$ is the projection,
     i.e.\ the following diagram is commutative.\smallskip
     \begin{center}
       \resetparms \qtriangle[X`Y\times \I`Y;g`f`P_Y]
     \end{center}
   \item[(b)] $g(H)=Y\times \I$ and each $g^{-1}(y,c)\cap H$, $(y,c)\in Y\times
     \mathfrak{C}$, is compact and $0$-dimensional.
   \end{enumerate}
   In particular, $H_{\mathfrak{C}}=H\cap g^{-1}(Y\times \mathfrak{C})$ is a
   closed section for $f$ such that $f\uhr H_{\mathfrak{C}}$ is a
   compact light map.
 \end{theorem}

 To prepare for the proof of Theorem \ref{theorem-cantor-bouquet}, we
 introduce some terminology. For a metric space $(X,d)$, a nonempty
 subset $A\subset X$ and $\varepsilon>0$, as in Section \ref{section-bula-property},
 we let
\[
B_\varepsilon^d(A)=\{x\in X: d(x,A)<\varepsilon\}.
\]
Also, we will use $\diam_d(A)$ to denote the \emph{diameter} of $A$
with respect to $d$.\medskip

Following \cite{gutev:93b}, to every nonempty subset $F\subset X$ we
associate the number
\[
\begin{split}
  \delta(F,X)=\inf\big\{1,\varepsilon: &\ \varepsilon>0\ \text{and}\\ &\ F\subset B_\varepsilon^d(S)\ \text{for some
  nonempty finite}\ S\subset F\big\}.
\end{split}
\]

In what follows, we let $\Omega(f)$ to be the set of all open sections for
$f$. Also, for every $U\in \Omega(f)$, we introduce the \emph{$d$-mesh} of
$U$ with respect to $f$ by letting $\mesh_d(U,f)=\sup\big\{\delta(f^{-1}(y)\cap U,
X):y\in Y\big\}$.

\begin{proposition}
  \label{proposition-disjoint-sections}
  Let $(X,d)$ be a metric space, $Y$ be a paracompact $C$-space, and
  let $f:X\to Y$ be an open continuous surjection such that each fiber
  of $f$ is dense in itself and $d$-complete. Then, for every $U\in
  \Omega(f)$ there are disjoint open sections $U_0,U_1\in \Omega(f)$ such that
  $\overline{U_i}\subset U$, $i=0,1$.
\end{proposition}

\begin{proof}
  Consider $U$ endowed with the compatible metric
  \[
  \rho(x,y)=d(x,y)+\left|\frac 1{d(x,X\setminus U)} - \frac 1{d(y,X\setminus
      U)}\right|,\quad x,y\in U.
  \]
  Next, define an l.s.c.\ mapping $\Phi:Y\smap X$ by $\Phi(y)=f^{-1}(y)\cap U$,
  $y\in Y$. Then, each $\Phi(y)$, $y\in Y$, is infinite and $\rho$-complete
  in $U$ because each $f^{-1}(y)$, $y\in Y$, is dense in itself and
  $d$-complete. Hence, by Corollary
  \ref{corollary-disjoint-michael-pair}, $\Phi$ has compact-valued
  u.s.c.\ multiselections $\psi_0,\psi_1:Y\smap U$ such that $\psi_0(y)\cap
  \psi_1(y)=\emptyset$ for every $y\in Y$. In fact, $\psi_0$ and $\psi_1$ are
  compact-valued and u.s.c.\ as mappings from $Y$ into the subsets of
  $X$. Hence, each $F_i=\bigcup\{\psi_i(y): y\in Y\}$, $i=0,1$, is a closed subset
  of $X$, with $F_i\subset U$ and $f(F_i)=Y$. Since $F_0\cap F_1=\emptyset$, we can
  take disjoint open sets $U_0, U_1\subset X$ such that $F_i\subset U_i\subset
  \overline{U_i}\subset U$, $i=0,1$. This completes the proof.
\end{proof}

In our next considerations, to every nonempty subset $F$ of a metric
space $(X,d)$ we associate (the possibly infinite) number
\[
\td_d(F)=\sup\big\{\diam_d(C): C\subset F\ \text{is connected}\big\}.
\]
Next, for a surjective map $f:X\to Y$ and a section $U\in \Omega(f)$, we let
\[
\td_d(U,f)=\sup\big\{\td_d(f^{-1}(y)\cap U): y\in Y\big\}.
\]
In the proof of our next lemma and in the sequel, $\omega$ denotes the
first infinite ordinal.

\begin{lemma}
  \label{lemma-small-sections}
  Let $(X,d)$ be a metric space, $Y$ be a paracompact $C$-space, and
  let $f:X\to Y$ be an open continuous surjection. Then, for every
  $\varepsilon>0$, every $G\in \Omega(f)$ contains an $U\in \Omega(f)$, with $\mesh_d(U,f)\leq \varepsilon$
  and $\td_d(U,f)\leq \varepsilon$.
\end{lemma}

\begin{proof}
  Let $\varepsilon>0$ and $G\in \Omega(f)$. Whenever $y\in Y$ and $n<\omega$, take an open
  subset $W_y^n\subset G$ such that $y\in f(W_y^n)$ and $\diam_d(W_y^n)<\varepsilon\cdot
  2^{-(n+1)}$. Since $f$ is open, each family
  $\mathscr{W}_n=\big\{f(W_y^n):y\in Y\big\}$, $n<\omega$, is an open cover of
  $Y$. Since $Y$ is a paracompact $C$-space, there now exists a
  sequence $\{\mathscr{V}_n:n<\omega\}$ of pairwise disjoint open families of
  $Y$ such that each $\mathscr{V}_n$, $n<\omega$, refines $\mathscr{W}_n$
  and $\mathscr{V}=\bigcup\{\mathscr{V}_n:n<\omega\}$ is a locally-finite cover of $Y$. For
  convenience, for every $n<\omega$, define a map $p_n:\mathscr{V}_n\to Y$ by
  $V\subset f\big(W^n_{p_n(V)}\big)$, $V\in \mathscr{V}_n$, and next set
  $U_{p_n(V)}=f^{-1}(V)\cap W^n_{p_n(V)}$.  We are going to show that
  \[
  U=\bigcup\big\{U_{p_n(V)}: V\in\mathscr{V}_n\ \text{and}\ n<\omega\big\}
  \]
  is as required. Since $\mathscr{V}$ is a cover of $Y$, $U$ is a
  section for $f$, and clearly it is open. Take a point $y\in Y$, and
  set $\mathscr{V}_y=\{V\in \mathscr{V}: y\in V\}$. Then, $\mathscr{V}_y$ is
  finite and $\big|\mathscr{V}_y\cap \mathscr{V}_n\big|\leq 1$ for every
  $n<\omega$ (recall that each family $\mathscr{V}_n$, $n<\omega$, is pairwise
  disjoint). Hence, we can numerate the elements of $\mathscr{V}_y$ as
  $\big\{V_k:k\in K(y)\big\}$ so that $V_k\in \mathscr{V}_k$, $k\in K(y)$,
  where $K(y)=\{n<\omega: \mathscr{V}_y\cap \mathscr{V}_n\neq\emptyset\}$. Next, set $U_k=
  U_{p_k(V_k)}$, $k\in K(y)$. Since
  \begin{equation}
    \label{eq:diameter-connected}
    \diam(U_k)<\varepsilon\cdot 2^{-(k+1)}\quad \text{for every $k\in K(y)$,}
  \end{equation}
  $f^{-1}(y)\cap U\subset B_\varepsilon^d(S)$ for every finite subset $S\subset f^{-1}(y)\cap U$,
  with $S\cap U_k\neq\emptyset$ for all $k\in K(y)$. Thus, $\delta(f^{-1}(y)\cap U,X)\leq \varepsilon$
  which completes the verification
  that $\mesh_d(U,f)\leq \varepsilon$.\smallskip

  To show that $\td_d(U,f)\leq \varepsilon$,
  take a nonempty connected subset $C\subset f^{-1}(y)\cap U$, and points
  $x,z\in C$. Since $C$ is connected and $C\subset \bigcup\{U_k:k\in K(y)\}$, there is a
  sequence $k_1,\dots, k_m$ of distinct elements of $K(y)$ such that
  $x\in U_{k_1}$, $z\in U_{k_m}$ and ${U_{k_i}\cap U_{k_j}\neq\emptyset}$ if and only if
  $|i-j|\leq 1$, see \cite[6.3.1]{engelking:89}. Therefore, by
  (\ref{eq:diameter-connected}),
  \begin{eqnarray*}
    d(x,z) &\leq & \sum_{i=1}^m\diam_d(U_{k_i})\leq \sum_{k\in K(y)}\diam_d(U_{k})\\
    &< &  \sum_{k\in K(y)}\varepsilon\cdot 2^{-(k+1)}< \varepsilon\cdot \sum_{k=0}^\infty 2^{-(k+1)}\quad =\varepsilon.
  \end{eqnarray*}
  Consequently, $\diam_d(C)\leq \varepsilon$, which completes the proof.
\end{proof}

Recall that a partially ordered set $(T,\preceq)$ is called a
\emph{tree} if the set $\{s\in T:s\prec t\}$ is well-ordered for every
point $t\in T$. Here, as usual, ``$s \prec t$'' means that $s\preceq
t$ and $s\neq t$. A \emph{chain} $\eta$ in a tree $(T,\preceq)$ is a
subset $\eta \subset T$ which is linearly ordered by $\preceq$. A maximal
chain $\eta$ in $T$ is called a \emph{branch} in $T$. Through this
paper, we will use $\mathscr{B}(T)$ to denote the set of all branches
in $T$. Following Nyikos \cite{nyikos:99}, for every $t\in T$, we let
\begin{equation}
  \label{eq:branch-base}
  U(t)=\{\beta\in \mathscr{B}(T): t\in\beta\},
\end{equation}
and next we set $\mathscr{U}(T)=\{U(t):t\in T\}$. It is well-known that
$\mathscr{U}(T)$ is a base for a non-Archimedean topology on
$\mathscr{B}(T)$, see \cite[Theorem 2.10]{nyikos:99}. In fact, one can
easily see that $s\prec t$ if and only if $U(t)\subset U(s)$, while $s$ and $t$
is incomparable if and only if $U(s)\cap U(t)=\emptyset$. In the sequel, we will
refer to $\mathscr{B}(T)$ as a \emph{branch space} if it is endowed
with this topology.\medskip

For a tree $(T,\preceq)$, let $T(0)$ be the set of all minimal elements of
$T$. Given an ordinal $\alpha$, if $T(\beta)$ is defined for every $\beta<\alpha$, then
we let
\[
T\uhr \alpha= \bigcup\{T(\beta):\beta\in \alpha\},
\]
and we will use $T(\alpha)$ to denote the minimal elements of $T\setminus
(T\uhr\alpha)$. The set $T(\alpha)$ is called the \emph{$\alpha^{\text{th}}$-level}
of $T$. The \emph{height} of $T$ is the least ordinal $\alpha$ such that
$T\uhr\alpha=T$. In particular, we will say that $T$ is an \emph{$\alpha$-tree}
if its height is $\alpha$. Finally, we can also define the \emph{height} of
an element $t\in T$, denoted by $\hgt(t)$, which is the unique ordinal
$\alpha$ such that $t\in T(\alpha)$.\medskip

In what follows, we will be mainly interested in $\omega$-trees, and the
following realization of the Cantor set as a branch space. Namely, let
$S$ be a set which has at least 2 distinct points, $S^n$ be the set of
all maps $t:n\to S$ (i.e., the $n^{\text{th}}$-power of $S$), and let
\[
S^{<\omega}=\bigcup\{S^{n+1}:n<\omega\}.
\]
Whenever $t\in S^{<\omega}$, let $\dom(t)$ be the \emph{domain} of
$t$. Consider the partial order $\preceq$ on $S^{<\omega}$ defined for $s,t\in
S^{<\omega}$ by $s\preceq t$ if and only if
\[
\dom(s)\subset \dom(t)\quad\text{and}\quad  t\uhr \dom(s)=s.
\]
Then, $(S^{<\omega},\preceq)$ is an $\omega$-tree such that its branch space
$\mathscr{B}(S^{<\omega})$ is the Baire space $S^\omega$. In particular, the
branch space $\mathscr{B}(2^{<\omega})$ is the Cantor set $2^\omega$. In the
sequel, we will refer to the tree $(2^{<\omega},\preceq)$ as the \emph{Cantor
  tree}.  \medskip

By Proposition \ref{proposition-disjoint-sections} and Lemma
\ref{lemma-small-sections}, using an induction on the levels of the
Cantor tree $(2^{<\omega},\preceq)$, we get the following immediate consequence.

\begin{corollary}
  \label{corollary-cantor-set-sections}
  Let $(2^{<\omega},\preceq)$ be the Cantor tree, and let $(X,d)$, $Y$, $f:X\to Y$
  and $U\in \Omega(f)$ be as in Theorem \ref{theorem-cantor-bouquet}. Then,
  there exists a map $h:2^{<\omega}\to \Omega(f)$ such that, for every two
  distinct members $s,t\in 2^{<\omega}$,
  \begin{enumerate}
  \item[(a)] $\overline{h(t)}\subset h(s)\subset U$ if $s\prec t$,
  \item[(b)] $h(s)\cap h(t)=\emptyset$ if $s$ and $t$ are incomparable,
  \item[(b)] $\mesh_d(h(t),f)\leq 2^{-\hgt(t)}$\ \ and\ \ $\td_d(h(t),f)\leq
    2^{-\hgt(t)}$.
  \end{enumerate}
\end{corollary}

We finalize the preparation for the proof of Theorem
\ref{theorem-cantor-bouquet} with the following special case of it.

\begin{lemma}
  \label{lemma-pactorization-cantor}
  Let $(X,d)$, $Y$, $f:X\to Y$ and $U\in \Omega(f)$ be as in Theorem
  \ref{theorem-cantor-bouquet}. Then, there exists a closed section
  $H\subset X$ of $f$, with $H\subset U$, and a surjective compact light map $\ell:H\to
  Y\times \mathfrak{C}$ such that $f\uhr H=P_Y\circ \ell$, i.e.\ the following
  diagram is commutative.\smallskip
  \begin{center}
    \resetparms \qtriangle[H`Y\times \mathfrak{C}`Y;\ell`f\uhr H`P_Y]
  \end{center}
  In particular, $f\uhr H$ is also a compact light map.
\end{lemma}

\begin{proof}
  Let $h:2^{<\omega}\to \Omega(f)$ be as in Corollary
  \ref{corollary-cantor-set-sections}. Whenever $n<\omega$, consider the
  $n^{\text{th}}$-level of the Cantor tree $2^{<\omega}$, which is, in
  fact, $2^{n+1}$. Then, set $H_n=h(2^{n+1})$, $n<\omega$, and
  $H=\bigcap\{H_n:n<\omega\}$. By (a) of Corollary
  \ref{corollary-cantor-set-sections}, $\overline{H_{n+1}}\subset H_n\subset U$
  for every $n<\omega$ because each level of $2^{<\omega}$ is finite. Hence, $H$
  is a closed subset of $X$, with $H\subset U$. Let us see that $H$ is a
  section for $f$. Indeed, take a point $y\in Y$, and a branch $\beta\in
  \mathscr{B}(2^{<\omega})$. Then, each $H_t(y)=h(t)\cap f^{-1}(y)$, $t\in \beta$,
  is a nonempty subset of $f^{-1}(y)$ (because $h(t)\in \Omega(f)$) such that
  $\overline{H_t(y)}\subset H_s(y)$ for $s\prec t$ (by (a) of Corollary
  \ref{corollary-cantor-set-sections}) and $\lim_{t\in\beta}\delta(H_t(y),X)=0$
  (by (c) of Corollary \ref{corollary-cantor-set-sections}). Hence, by
  \cite[Lemma 3.2]{gutev:93b}, $H_\beta(y)=\bigcap\{H_t(y): t\in \beta\}$ is a nonempty
  compact subset of $X$. Clearly, $H_\beta(y)\subset H\cap f^{-1}(y)$ which
  completes the verification that $H$ is a section for $f$. In fact,
  this defines a compact-valued mapping $\varphi:Y\times \mathscr{B}(2^{<\omega})\smap
  H$ by letting $\varphi(y,\beta)=H_\beta(y)=\bigcap\{h(t)\cap f^{-1}(y): t\in \beta\}$, $(y,\beta)\in Y\times
  \mathscr{B}(2^{<\omega})$. Since $\varphi(y,\beta)\subset f^{-1}(y)$, $(y,\beta)\in Y\times
  \mathscr{B}(2^{<\omega})$, the mapping $\varphi$ is the inverse $\ell^{-1}$ of a
  surjective single-valued map $\ell:H\to Y\times \mathscr{B}(2^{<\omega})$. Also,
  $\ell(x)=(y,\beta)$ if and only if $x\in \varphi(y,\beta)\subset f^{-1}(y)$, hence $f\uhr
  H=P_Y\circ \ell$.\medskip

  To show that $\ell$ is continuous and light, take an open set $V\subset Y$,
  $t\in 2^{<\omega}$, and let $U(t)$ be as in \eqref{eq:branch-base}. Then,
  $h(t)$ is an open set in $X$ such that, by (b) of Corollary
  \ref{corollary-cantor-set-sections}, $\ell^{-1}(y,\beta)=\varphi(y,\beta)\subset h(t)$ if
  and only if $t\in \beta$ (i.e., $\beta\in U(t)$). Consequently, $\ell^{-1}(V\times
  U(t))=f^{-1}(V)\cap h(t)\cap H$ is open in $H$. Finally, take a nonempty
  connected subset $C\subset \ell^{-1}(y,\beta)= \varphi(y,\beta)$ for a point $y\in Y$ and a
  branch $\beta\in \mathscr{B}(2^{<\omega})$. Then, $C\subset h(t)\cap f^{-1}(y)$ for
  every $t\in\beta$ and therefore, by (c) of Corollary
  \ref{corollary-cantor-set-sections}, $\diam_d(C)=0$. Hence, $C$ is a
  singleton, which implies that $\ell^{-1}(y,\beta)$ is $0$-dimensional
  being compact.\medskip

  To show finally that $f\uhr H$ is a compact light map, take a point
  $y\in Y$, and let us observe that $\ell\uhr \big(f^{-1}(y)\cap H\big)$ is
  perfect. Indeed, take a branch $\beta\in \mathscr{B}(2^{<\omega})$ and a
  neighbourhood $W$ of $\ell^{-1}(y,\beta)$ in $X$. Then, by \cite[Lemma
  3.2]{gutev:93b}, there exists a $t\in\beta$, with $H_t(y)=h(t)\cap f^{-1}(y)\subset
  W$. In this case, $\ell^{-1}(y,\gamma)\subset W$ for every $\gamma\in U(t)$, where $U(t)$
  is as in (\ref{eq:branch-base}). Namely, $\gamma\in U(t)$ implies that
  $t\in\gamma$ and, therefore, $\ell^{-1}(y,\gamma)\subset H_t(y)\subset W$. Thus, $\ell\uhr
  \big(f^{-1}(y)\cap H\big)$ is perfect, which implies that $f^{-1}(y)\cap
  H=\ell^{-1}(\{y\}\times \mathscr{B}(2^{<\omega}))$ is compact because so is
  $\mathscr{B}(2^{<\omega})$. Since $\mathscr{B}(2^{<\omega})$ is
  zero-diemnsional and $\ell$ is a light map, according to the classical
  Hurewicz theorem (see, \cite{engelking:95}), this also implies that
  $\dim(f^{-1}(y)\cap H)=0$ which completes the proof.
\end{proof}

\begin{proof}[Proof of Theorem \ref{theorem-cantor-bouquet}]
  We repeat the arguments of \cite[Theorem 1]{bula:83}. Briefly, let
  $(X,d)$, $Y$, $f:X\to Y$ and $U\in \Omega(f)$ be as in Theorem
  \ref{theorem-cantor-bouquet}. By Lemma
  \ref{lemma-pactorization-cantor}, there exists a closed section $H\subset
  X$ of $f$, with $H\subset U$, and a continuous surjective map $\ell:H\to Y\times
  \mathfrak{C}$ such that $f\uhr H$ is a compact light map, and $f\uhr
  H=P_Y\circ \ell$. Take a continuous surjective map $\rho:\mathfrak{C}\to \I$
  such that the set
  \[
  D=\{t\in \I: |\rho^{-1}(t)|>1\}
  \]
  is countable. Also, let $P_{\mathfrak{C}}:Y\times \mathfrak{C}\to
  \mathfrak{C}$ be the projection. Then, using the Tietze-Urysohn
  theorem, extend $\rho\circ P_{\mathfrak{C}}\circ \ell$ to a continuous map $u:X\to
  \I$. In this way, we have that
  \[
  u(f^{-1}(y)\cap H)=\I\quad \text{for every $y\in Y$.}
  \]
  Then, we can define our $g:X\to Y\times \I$ by $g(x)=(f(x),u(x))$, $x\in
  X$. As for the the second part of Theorem
  \ref{theorem-cantor-bouquet}, take a copy $\mathfrak{C}$ of the
  Cantor set in $\I\setminus D$, which is possible because $D$ is
  countable. Then, by the properties of $\rho$, we have that $g^{-1}(Y\times
  \{c\})\cap H= \ell^{-1}(Y\times\{c\})$ for every $c\in \mathfrak{C}$. Hence, Lemma
  \ref{lemma-pactorization-cantor} completes the proof.
\end{proof}

We finalize this paper with several applications of Theorem
\ref{theorem-cantor-bouquet}. In what follows, for a space $X$, let
$\mathscr{F}(X)$ be the set of all nonempty closed subsets of
$X$. Recall that the \emph{Vietoris topology} $\tau_V$ on
$\mathscr{F}(X)$ is generated by all collections of the form
$$
\langle\mathscr{V}\rangle = \left\{S\in \mathscr{F}(X) : S\subset
  \bigcup \mathscr{V}\ \ \text{and}\ \ S\cap V\neq \emptyset,\
  \hbox{whenever}\ V\in \mathscr{V}\right\},
$$
where $\mathscr{V}$ runs over the finite families of open subsets of
$X$. In the sequel, any subset $\mathscr{D}\subset \mathscr{F}(X)$ will
carry the relative Vietoris topology $\tau_V$ as a subspace of
$(\mathscr{F}(X),\tau_V)$. In fact, we will be mainly interested in the
subset
\[
\mathscr{F}(f)=\{H\in \mathscr{F}(X): f(H)=Y\},
\]
where $f:X\to Y$ is a surjective map.

\begin{corollary}
  \label{corollary-closed-section-light-dense}
  Let $(X,d)$ be a metric space, $Y$ be a paracompact $C$-space, and
  let $f:X\to Y$ be an open continuous surjection such that each fiber
  of $f$ is dense in itself and $d$-complete. Then, the set
  \[
  \mathscr{L}(f)=\big\{H\in \mathscr{F}(f): f\uhr H\ \text{is a compact
    light map}\big\}
  \]
  is dense in $\mathscr{F}(f)$ with respect to the Vietoris topology
  $\tau_V$.
\end{corollary}

\begin{proof}
  Take a closed section $F\in \mathscr{F}(f)$, and a finite family
  $\mathscr{U}$ of open subsets of $X$, with $F\in  \langle
  \mathscr{U}\rangle$. Then, $U=\bigcup\mathscr{U}$ is an open section for $f$,
  so, by Theorem \ref{theorem-cantor-bouquet}, it contains a closed
  section $H\subset U$ such that $f\uhr H$ is a compact light map. Take
  a finite set $S\in\langle\mathscr{U}\rangle$, and then set $Z=H\cup S$. Clearly,
  $Z\in \mathscr{L}(f)\cap \langle\mathscr{U}\rangle$, which completes the proof.
\end{proof}

\begin{proposition}
  \label{proposition-zero-set-cantor}
  Whenever $Y$ is a metrizable space, there exists a closed
  $0$-dimensional subset $A\subset Y\times \mathfrak{C}$ such that $P_Y(A)=Y$.
\end{proposition}

\begin{proof}
  We follow the idea of \cite[Lemma 4.1]{tuncali-valov:02}. Fix a
  $0$-dimensional metrizable space $M$ and a perfect
  surjective map $h\colon M\to Y$.  By
  \cite[Proposition 9.1]{pasynkov:98}, there exists a continuous map
  $g\colon M\to Q$, where $Q$ is the Hilbert cube, such that
  $h\triangle g\colon M\to Y\times Q$ is an
  embedding. Next, take a Milyutin map $p\colon \mathfrak{C}\to Q$,
  i.e.\ a surjective continuous map admitting an averaging operator
  between the function spaces $C(\mathfrak{C})$ and $C(Q)$, see
  \cite{pelczynski:68}.  According to \cite{ditor:73}, there exists a
  compact-valued lower semi-continuous map $\varphi\colon Q\smap
  \mathfrak{C}$ such that $\varphi (z)\subset p^{-1}(z)$ for all $z\in
  Q$. Applying Michael's $0$-dimensional selection theorem
  \cite{michael:56}, there exists a continuous
  map $\ell:M\to \mathfrak{C}$, with $\ell(x)\in\varphi(g(x))$ for any
  $x\in M$. Then, $h\triangle \ell$ embeds $M$ as a closed subset $A$ of
  $Y\times \mathfrak{C}$. Obviously,
  $A$ is $0$-dimensional and
  $P_Y(A)=Y$.
\end{proof}

\begin{corollary}
  \label{corollary-closed-section-zero-dense}
  Let $X$ be a metrizable space, $Y$ be a metrizable $C$-space, and
  let $f:X\to Y$ be an open perfect surjection such that each fiber
  of $f$ is dense in itself. Then, the set
  \[
  \mathscr{F}_0(f)=\big\{H\in \mathscr{F}(f): \dim(H)=0\big\}
  \]
  is dense in $\mathscr{F}(f)$ with respect to the Vietoris topology
  $\tau_V$.
\end{corollary}

\begin{proof}
  Take a closed section $F\in \mathscr{F}(f)$, and a finite family
  $\mathscr{U}$ of open subsets of $X$, with $F\in \langle
  \mathscr{U}\rangle$. Then, $U=\bigcup\mathscr{U}$ is an open section for $f$,
  so, by Theorem \ref{theorem-cantor-bouquet}, there exists a closed
  section $H$ for $f$, with $H\subset U$, a continuous surjective map $g:X\to
  Y\times \I$, and a copy $\mathfrak{C}\subset \I$ of the Cantor set such that
  $f=P_Y\circ g$, $g(H)=Y\times \I$, and $f\uhr \big(H\cap g^{-1}(Y\times
  \mathfrak{C}\big)$ is a light map. By Proposition
  \ref{proposition-zero-set-cantor}, $Y\times \mathfrak{C}$ contains a
  closed $0$-dimensional set $A$, with $P_Y(A)=Y$. Finally,
  take $B=H\cap g^{-1}(A)$ which is a closed section for $f$ because
  $P_Y(A)=Y$. Since $f$ is perfect, so is $g$. Hence, $g\uhr B$ is a
  perfect light map and, according to the classical
  Hurewicz theorem, $\dim(B)=0$. Then, $Z=B\cup S\in \mathscr{F}_0(f)\cap \langle
  \mathscr{U}\rangle $ for some (every) finite set $S\in\langle\mathscr{U}\rangle$.
\end{proof}

To prepare for our last consequence, let us also observe the following
property of $0$-dimensional sections.

\begin{proposition}
  \label{proposition-g-delta-zero-dense}
  Let $X$ be a compact metrizable space, and let $\mathscr{F}_0(X)$ be
  the subset of all $0$-dimensional members of $\mathscr{F}(X)$. Then,
  $\mathscr{F}_0(X)$ is a $G_\delta$-subset of $\mathscr{F}(X)$.
\end{proposition}

\begin{proof}
  Take a metric $d$ on $X$ compatible with the topology of $X$. Next,
  for every $H\in \mathscr{F}_0(X)$ and $n\geq 1$, take a pairwise disjoint
  finite family $\mathscr{V}_n(H)$ of open subsets of $X$ such that
  $\diam(V)<1/n$, $V\in \mathscr{V}_n(H)$, and $H\in \langle
  \mathscr{V}_n(H)\rangle$. Then, each
  $\mathscr{V}_n=\bigcup\big\{\mathscr{V}_n(H): H\in \mathscr{F}_0(X)\big\}$,
  $n\geq 1$, is $\tau_V$-neighbourhood of $\mathscr{F}_0(X)$, and clearly
  $\mathscr{F}_0(X)= \bigcap\big\{\mathscr{V}_n(F): n=1,2,\dots\big\}$.
\end{proof}

According to Corollary \ref{corollary-closed-section-zero-dense} and
Proposition \ref{proposition-g-delta-zero-dense}, we have the
following immediate consequence, which is the Levin-Rogers
\cite[Theorem 1.2]{levin-rogers:99}.

\begin{corollary}
  [\cite{levin-rogers:99}]
  \label{corollary-closed-section-zero-compact}
  Let $X$ be a compact metrizable space, $Y$ be a metrizable
  $C$-space, and let $f:X\to Y$ be an open continuous surjection such
  that each fiber of $f$ is dense in itself. Then, the set
  \[
  \mathscr{F}_0(f)=\big\{H\in \mathscr{F}(f): \dim(H)=0\big\}
  \]
  is a dense $G_\delta$-subset in $\mathscr{F}(f)$ with respect to the
  Vietoris topology $\tau_V$.
\end{corollary}

\newcommand{\noopsort}[1]{} \newcommand{\singleletter}[1]{#1}
\providecommand{\bysame}{\leavevmode\hbox to3em{\hrulefill}\thinspace}
\providecommand{\MR}{\relax\ifhmode\unskip\space\fi MR }
\providecommand{\MRhref}[2]{%
  \href{http://www.ams.org/mathscinet-getitem?mr=#1}{#2}
}
\providecommand{\href}[2]{#2}


\end{document}